\documentclass[12pt,oneside]{article}
\usepackage{verbatim}
\usepackage{graphicx}
\usepackage{amsmath}
\usepackage{amsfonts}
\usepackage{algorithmic}
\usepackage{algorithm}
\usepackage{rotating}
\usepackage{multicol}
\usepackage{enumerate}

\usepackage{pdflscape}
\usepackage{xcolor}
\usepackage[normalem]{ulem}

\usepackage{natbib}
\usepackage{booktabs}
\usepackage[T1]{fontenc}
\usepackage{graphics}
\usepackage{amssymb}
\usepackage{pgf}
\usepackage{colortbl}
\usepackage{multirow}
\usepackage{booktabs}

\newcommand{\red}[1]{\textcolor{black}{#1}}

\begin{document}

\title{A simulation framework for a station-based bike-sharing system}
\author{\it E. Angelelli$^1$, M. Chiari$^2$, A. Mor$^3$, M.G. Speranza$^1$ \\
{\small \it $^1$Department of Economics and Management}\\
{\small \it University of Brescia, Italy}\\
{\small \{enrico.angelelli, grazia.speranza\}@unibs.it}\\
{\small \it $^2$Brescia Mobilità S.p.A.}\\
{\small \it Brescia, Italy}\\
{\small mchiari@bresciamobilita.it}\\
{\small \it $^3$Department of Management, Economics and Industrial Engineering}\\
{\small \it Politecnico di Milano, Italy}\\
{\small andrea.mor@polimi.it}\\
}
\date{}
\maketitle

\begin{abstract}
Many cities and towns offer to citizens a Bike-Sharing System (BSS). When a company starts the service, multiple decisions have to be taken \red{at strategic, tactical, and operational level. When the service is in place, it is often necessary to adapt these decisions to changed conditions.
In this paper, starting from the experience gained in the real-case of Bicimia in Brescia, Italy, we present a simulation framework to support decisions in the design or revision of a BSS, including the shifts of the operators, the number of vehicles and their capacity, the sizing of the stations, and the sizing of the fleet of bikes}.
\vspace{3mm}

\end{abstract}

\section{Introduction}\label{sec:intro}

The first Bike-Sharing Systems (BSS) date back more than 50 years and have seen a great boost in their popularity in the last 15 years, becoming increasingly important in urban environments as a low-cost, environmentally-friendly transportation means for short-distance journeys (see \cite{demaio2009bike} for a historical overview of BSSs).
Two types of BSS have emerged over the years: station-based (or docked) BSS and free-floating (or dockless) BSS.
In station-based BSSs, bike stations are scattered in the area served by the system. Each station is equipped with a number of stands, each holding a bike or vacant. Bikes can be rented at a station and returned to a different station.
In free-floating BSSs, no stations are present and bikes can be rented or returned in any point of the area in which the system operates. While free-floating BSSs are gaining popularity, the majority of BSSs currently operating around the world are station-based (see \cite{meddin2021bike} and \cite{wiki}). For this reason, together with the fact that this work has its foundation in a collaboration with a company operating a station-based BSS, in the following we focus on this type of systems.

When a company starts a station-based BSS, the decisions to be taken include the stations location and the number of stands for each of them, the number of bikes in the system, the number of vehicles for the rebalancing operations to mitigate imbalances in bike requests (i.e., transferring bikes from stations with too many bikes to stations with too few), and the shifts of the operators, in addition to the cost and the method for the payment of the service. When the service is in place, it is often necessary to adapt it to changed conditions and several decisions have to be revised.
The planning problems arising in a BSS at the strategic, tactical, and operational level have been discussed in \cite{vogel2016service}, \cite{espegren2016static}, and recently in \cite{shui2020review}.
A study on the strategic design of BSSs with service level constraints is considered in \cite{lin2011strategic} to determine the number and the location of bike stations and the network structure of the paths connecting the stations. The planning activities arising at the operational level include bike relocation. Challenges and opportunities in station-based bike-sharing relocation are presented in \cite{vallez2021challenges}.

In this paper, we present a simulation framework to support decision making regarding bike station sizing, bike fleet sizing, vehicle fleet sizing and rostering for a station-based BSS where the bike relocation operations are performed dynamically during the day.
The framework simulates a day of service and requires, as input data, the location of the stations and, for each station, the number of stands and the number of bikes at the beginning of the day. It also requires the capacity of each vehicle in the fleet and the shifts of the operators. Requests of bikes and vacant stands in each station are modeled, starting from real data, as non-homogeneous Poisson processes.
In general, strategic and tactical decisions, such as the number of vehicles, the sizing of the stations, and the work shifts of the operators, are strongly connected to the operational ones, that is the efficiency of the routes of the vehicles performing the bike relocation and the number of bikes they transfer among the stations.
If the operations are optimized, a better service level can be achieved with the available resources or less resources are needed to achieve a given service level.
For this reason, our framework simulates a day of service while optimizing the operations of the vehicles with the goal of minimizing the shortages (i.e., lost rental and return requests due to the lack of bikes and vacant stands, respectively). Within the presented simulation framework, bike relocation is performed dynamically during the day. The relocation problem is modeled as a Dynamic Bike Relocation Problem (DBRP) (see \cite{shui2020review} for references).
A reoptimization approach is presented for this problem, where reoptimization epochs are defined dynamically based on the actions of the vehicles, which are guided by the status of the system.
A look-ahead algorithm is adopted for the reoptimization problem, using the forecast of the requests of bikes and vacant stands to anticipate the consequences of the operations and reduce the shortages.

The output of the simulation consists of indicators of the performance of the BSS.
As the information about unmet rentals and returns cannot be directly measured by the system, 
one of the key performance metrics is the total amount of time stations are either full or empty. This metric represents an upper bound on the shortage time of the system as, for instance, a bike shortage can only happen if a station is empty.
Our simulation framework allows the company in charge of the service to explore 'what if' questions and scenarios without having to experiment on the system itself and find the trade-off between the cost of the service and its quality, as perceived by the users. 

The simulation framework we present has been developed during a collaboration with Brescia Mobilità, the public  company that is responsible for several mobility services in Brescia, Italy, including the station-based BSS called Bicimia. The experiments of the simulation framework we present in this paper are based on real data from Bicimia.  In general, and also in the case of Bicimia, the level of usage of a BSS during a day strongly depends on the month (January is very different from August), the day of the week (working or festive day), and  on the weather conditions (rainy or sunny day) (see \cite{eren2020review}).
In the run experiments the forecast of the requests of bikes and vacant stands in each station is obtained by using historical data. 

	The topic of simulating a BSS has been tackled in various ways in the literature.
	\cite{romero2012simulation} propose a simulation tool integrating private cars and shared bikes. The tool is used to optimize the location of stations for the BSS in order to improve the efficiency and sustainability of the transport system.  
	\cite{jian2016simulation} present a simulation optimization approach for large-scale BSSs to set the capacity of each station and the number of bikes available in each station at the beginning of the day with the aim of minimizing the expected number of unsatisfied customers.

    Below are the contributions that are particularly relevant to the work presented in this paper.
	
	In \cite{soriguera2018simulation}  an agent-based simulation framework is proposed. Both periodic and continuous relocation algorithms are tested in the simulator. Contrary to our contribution, the proposed algorithms consider only present information and do not have any look-ahead capability.
	An agent-based simulation environment for a station-based BSS is also presented in \cite{fernandez2020bike3s}. The system grants the users the possibility of reservation. The user demand is generated through a homogeneous Poisson distribution. The tool is validated on data from BSSs in Madrid and London. While the tool is described as capable of considering relocation policies, no specific detail is provided.
	\cite{jin2022simulation} present a simulation framework for a station-based BSS considering valet service for high-volume stations at peak hours and bike breakdown. Both static and dynamic relocation algorithm are tested. This contribution differs from the present one in that the only dynamic rebalancing approach proposed has fixed reoptimization epochs, whereas the presented approach considers reoptimization epochs that are dependent on the behavior of the vehicles (which itself is guided by the status of the system). Furthermore, in \cite{jin2022simulation} vehicle interaction is not considered when performing relocation: stations are clustered and each cluster is assigned to a vehicle. Finally, the routing problem faced when relocating bikes minimizes the sum of the travel time and loading/unloading operations time, as opposed to the minimization of bike and stand shortages proposed in this contribution.

	Recent research on the DBRP has been highlighted in \cite{shui2020review}. In the following, we report the contributions that are particularly relevant for this paper.
	In \cite{chemla2013self} the DBRP for a station-based BSS is investigated. User requests are assumed to follow a homogeneous Poisson process defined for each station. The destination of each user request is also drawn at random, according to a given probability matrix.
	\cite{caggiani2013dynamic} present a simulation model for the dynamic relocation of bikes to minimize the relocation costs while ensuring a high level of user satisfaction. Days are assumed to be divided into discrete intervals, each with different user demand. Different heuristics for the relocation are tested using the discrete-events simulation framework.
	A study on a proactive approach to demand fluctuation in the context of a DBRP is presented in \cite{regue2014proactive}. A forecasting model is used to anticipate the inventory level for different time horizons (20, 40, and 60 minutes). Forecasts are used together with an inventory level model to define the redistribution needs. The routing component of the problem is tackled by preprocessing the set of stations to be visited by each vehicle and solving a single-vehicle routing problem for each.
	In \cite{legros2019dynamic} a dynamic strategy for the DBRP is studied. The routing problem is modeled through a Markov Decision Process. Similarly to our contribution, the arrival rate of user rentals and returns is modeled by a non-homogeneous Poisson process. However, drawing from the case of Velib' in Paris, each relocating vehicle is assigned to an area and vehicles do not interact. In \cite{brinkmann2019dynamic}  a dynamic lookahead policy for the stochastic dynamic inventory routing problem is presented to minimize the expected amount of unsatisfied demand. One vehicle is used to perform the relocation of bikes.
	\cite{jimenez2020new} expand the work of \cite{soriguera2018simulation} discussing a three stage algorithm. First, the ideal bike distribution is calculated for each station. Then, optimal repositioning routes are scheduled and, finally, these routes are modified according to minimum-cost criteria. Bike returns and rentals at each station are treated as independent Poisson events and the accumulated number of rentals and returns is approximated by a Normal distribution.

It is also worth mentioning that various works on the data available about the behavior of users in a BSS have been published in recent years. We report those of particular interest and relevance to this paper. \cite{negahban2019simulation} presents a novel approach to overcome the intrinsic censoring of user demand and estimate the real demand starting from the partial information available. \cite{eren2020review} present a review of the factors influencing user demand, in particular weather conditions and temporal factors. Various studies have also been conducted to model user behavior and spatio-temporal patterns in the usage of a BSS, as for instance in \cite{corcoran2014spatio}, \cite{zhang2018mining}, \cite{kutela2019influence}, and \cite{ashqar2019modeling}.

\bigskip

\red{The contributions of this paper can be summarized as follows:}
\begin{itemize}
	\item \red{a simulation framework for a BSS is introduced which accounts for its setting (e.g., the fleet size, the location of stations), as well as the users' requests
	and the forecast of such requests according to the characteristics of the day being simulated;}
	\item \red{daily operations for bike relocation are modeled as a DBRP. An algorithm is proposed for the solution of the problem where
	reoptimization epochs are defined dynamically based on the actions of the vehicles, which operate in the same area.
	The forecast of the users' requests in the system is used to anticipate the consequences of the operations and reduce shortages;}
	\item \red{a real case is presented and results provided by application of the framework are discussed in terms of the shifts of the operators, the number of vehicles and their capacity, the sizing of the stations, and the sizing of the fleet of bikes.}
\end{itemize}

The paper is structured as follows. In Section \ref{sec:tool} the  proposed simulation framework is presented. The case study is discussed in Section \ref{sec:casestudy} and computational results are presented in Section \ref{sec:compres}. Finally, conclusions are drawn in Section  \ref{sec:conclusions}.

\section{Simulation framework}\label{sec:tool}

When a station-based BSS has to be started by a company or is already in place and has to be revised, the management has to take decisions about the sizing of the stations and that of the fleet of bikes, as well as the sizing of the fleet of repositioning vehicles and shift schedules. However, these strategic and tactical decisions are connected to the operational decisions of the fleet, like those regarding its routing and the inventory management of the stations.
With the aim to assist the strategic and tactical decision-making of the company, a simulation framework for a station-based BSS is presented in this section.
The simulation framework heavily depends on a representation of the 'on the field operations' which we will describe in details below.
First, in Section \ref{sec:framework} the structure of the framework is presented.
Then, in Section \ref{sec:algo} the optimization approach for the solution of the bike relocation problem embedded into the simulation framework is presented.

\subsection{Structure of the simulation framework} \label{sec:framework}

\begin{figure}[H]
	\centering
	\includegraphics[width=0.95\linewidth]{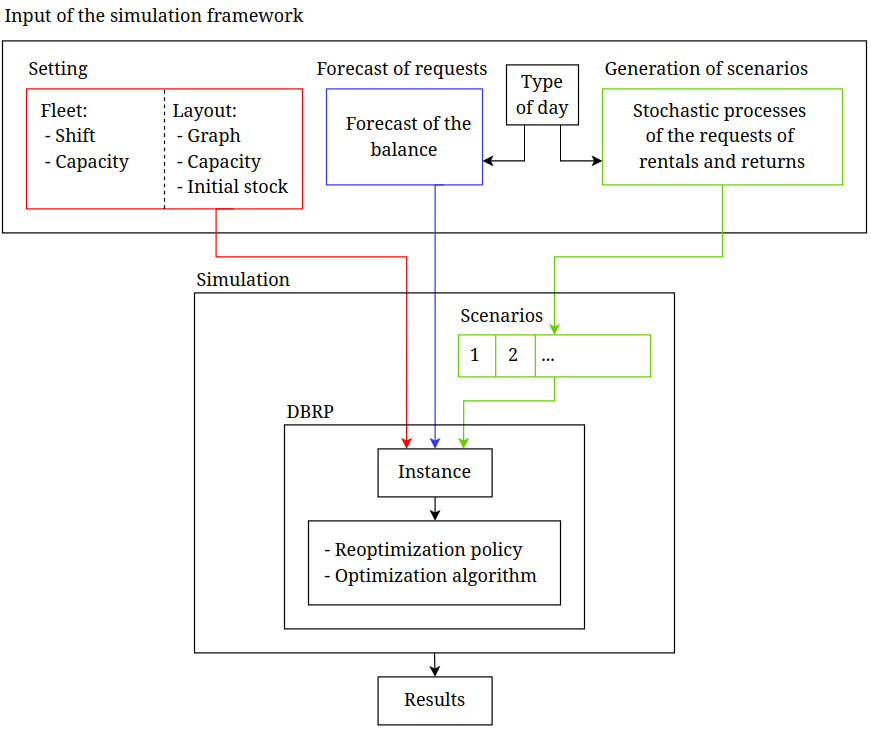}
	\caption{The structure of the simulation framework.}
	\label{fig:simulation}
\end{figure}

The structure of the simulation framework is illustrated in Figure \ref{fig:simulation}. To simulate the BSS, the framework requires several inputs. The first set of inputs is the setting of the system, that is the layout of the stations and details about the fleet of vehicles; the second set concerns the forecast of the requests of the users and is composed of the forecast of the balance (difference) between potential return and rental requests at each station; the third set is about the stochastic processes representing users requests. The sampling of such processes yields a set of scenarios. The forecast of requests and the generation of scenarios may be influenced by the type of day under consideration.
All these inputs are discussed in the remainder of this section.

\paragraph{Setting of the system.} We define the setting of the BSS by two major components: layout and fleet, which describe the operational environment for the simulation of the bike relocation policy. 
\subparagraph{Layout.}
We call \textit{layout} the data that describes the stations, their locations, the number of stands available in each station, and the initial stock of each station.
Formally, we define a layout as  $L = (G, U, S)$, where $G$ represents the road network connecting the stations, $U$ represents the stations capacity and $S$ represents the station stock level at the beginning of the day.
More in particular, $G$ is a directed graph $(V,A)$ with $V = \{0\} \cup V'$, where the vertex $0$ represents a depot, used to store the spare and broken bikes and as a parking space for the vehicles, and the vertices $V'=\{1,\ldots,n\}$  represent stations of the BSS.
An arc $(i,j)\in A$ exists between any pair of vertices, with $c_{ij}$ indicating the traveling time needed for a vehicle to travel directly from $i$ to $j$ and $d_{ij}$ the distance.
Each station $i$ is equipped with $u_i$ stands. Each stand can hold one bike. A stand can be used to make a bike available for rental or can be left vacant to make it possible for a user to return a bike. At any moment a station $i$ can offer at most $u_i$ bikes for rental if all the stands are taken or at most $u_i$ for the return of bikes if no stand is occupied by a bike. Thus, $U = \{u_i | i \in V'\}$.
Finally, the information system of the company collects real-time information on the actions performed by the users, i.e., rentals and returns. Real-time information on the number of bikes in each station is also available. We call this the \textit{stock} of a station. Within a day of operations, time is defined as continuous. Let $t \in [0, H]$ indicate a time instant, with $0$ being the beginning instant and $H$ being the final instant of the day. We denote the stock of station $i$ at time $t$ as $s_i(t)$. At time $t'$, information on $s_i(t)$ is available for $t \leq t'$, that is, the present and past stock of the station. The information system is also aware of the actions performed by the vehicles, that is, traveling, arriving at a station, loading or unloading bikes, and the departure from a station. Thus, $S = \{s_i(0) | i \in V'\}$.

\subparagraph{Fleet.}
We call \textit{fleet} the data about the vehicles, that is, the number of vehicles and, for each, its capacity and shift times.
Formally, we define a fleet as $M = \{(Q_k,\alpha_k,\beta_k) | k \in K\}$, 
where $K$ is the set of heterogeneous vehicles available to relocate bikes between stations during the day.
A vehicle $k \in K$ becomes available at the depot at the beginning of its shift $\alpha_k$, possibly loads some bikes, travels to a station, loads or unloads bikes, travels to a new station, repeats the sequence of operations, and returns to the depot at the end of its shift $\beta_k$.
Vehicle $k \in K$ has capacity $Q_k$ which means that the vehicle can transport $Q_k$ bikes at a time at most.

\paragraph{Type of day.} 

It is worth mentioning that the framework can consider different types of days characterized by a set of features like the type of weekday (i.e., working day, Saturday, Sunday), weather (sunny, rainy), and month (from January to December).

It is easy to see that the type of day may have a direct influence on the inclination of users to use the system. Thus, both the forecast of the balance and the stochastic processes representing the request are determined after the choice of a type of day.

The efficiency of a setting may also depend on the type of day, but setting choice is up to the management, who can plan different settings (e.g., shifts) according to the season or weekday.
According to data dependency on the features of the type of day, the simulator covers events happening in a single day and time spans in the interval $[0,H]$ where 0 corresponds to 00:00:00 and H to 24:00:00.

\paragraph{Forecast of the requests.} 

To take decisions at the operational level, such as the relocation of bikes among the stations, it is essential to be able to predict the behavior of the users of the BSS. Let us denote by $\bar{b}_{i}^+(t)$ and $\bar{b}_{i}^-(t)$ the forecast of the number of bikes that the users wish to return or rent, respectively, at station $i$ in the time interval $[0,t]$.
The forecast of the \textit{balance} between desired rentals and returns is thus defined as $\bar{f}_{i}(t) = \bar{b}_{i}^+(t) - \bar{b}_{i}^-(t)$. It must be specified that $\bar{f}_{i}(t)$ represents the difference of bikes that users wish to return and wish to rent in the interval $[0,t]$ regardless of the actual status of the station at time $t$, meaning that it might take negative values or values greater than the capacity of station $i$. This is because $\bar{f}_{i}(t)$ represents the need of bikes and does not take into account the stock and the capacity of a station. 

\paragraph{Generation of the scenarios.} 

We model the requests of the users by means of stochastic processes that describe the number of requests of rentals from and returns to each station during the day. Let $X_{i}^+(t)$ and $X_{i}^-(t)$ be stochastic processes that represent the number of requests of rentals from and returns to station $i$ up to time $t$, respectively.
We call {\em scenario} a realization $\{(x_i^+(t), x_i^-(t))| i\in V'\}$ of the stochastic processes $X_{i}^+(t)$ and $X_{i}^-(t)$, respectively, generated by sampling. Processes $X_{i}^+(t)$ and $X_{i}^-(t)$ have been assumed to be non-homogeneous Poisson processes. That is, for each station $i$, the rate of returns and rentals is described by rate $\lambda_{i}^+(t)$ and $\lambda_{i}^-(t)$, respectively, which we will assume to be non-negative stepwise constant functions. This allows to represent the variation of the rate of requests during the day.

\bigskip
For each scenario, the operational decisions taken by the system (and actions to be implemented by the fleet) are determined by solving a Dynamic Bike Relocation Problem (DBRP), which is discussed in Section \ref{sec:algo}. The performance of the system is reported as the average of the results over the set of scenarios. The decision-maker can use the simulation framework to test the performance of the service with different settings and to evaluate their impact on the performance of the service.

\subsection{Dynamic Bike Relocation Problem} \label{sec:algo}

In this section, the DBRP is defined to model the operational decisions of the company within the simulation framework. The policy defining the reoptimization epochs, the information available at each epoch, and the decisions to be taken in each epoch are first described. Then, the optimization algorithm used for the relocation operations is presented.

\paragraph{Reoptimization policy.}

The simulation framework uses a dynamic reoptimization strategy. Reoptimization epochs are the following:

\begin{itemize}
	\item the beginning of the first shift of the day being simulated, to build an initial plan for all the scheduled vehicles;
	\item every time a vehicle enters a station, to reevaluate its intervention based on the actual status of the station;
	\item every time a vehicle exits a station, to reevaluate the plan of the vehicles with the changes in the status of the stations that occurred while the operator was loading/unloading bikes (i.e., bikes were rented or returned in other stations of the system). As result, the vehicle will either stay at the station to perform further actions or depart for a different station.
\end{itemize}

It is worth stating that, while the two latter reoptimization epochs are triggered by the actions of one vehicle, the optimization algorithm builds a solution for all the vehicles whose shift has not already ended, based on the information available at the reoptimization epoch.

\paragraph{Information available at each epoch.}

In each epoch $t'$, the information available to plan for the relocation of bikes is the following. In addition to the setting of the system, the stock $s_i(t')$ is known for each station $i$.
In addition, information on the current load of the vehicles and on the operations being performed is also available. In particular, at any time $t$, the number of bikes on board of each vehicle is known, as well as if the vehicle is in a station performing loading/unloading operations, or if it is en route, traveling to a station.

\paragraph{Decisions to be taken.}

In each reoptimization epoch, the decisions to be taken concern the routing of the vehicles and the number of bikes to be loaded or unloaded at each station (or at the depot). We call the loading/unloading action of a vehicle an \textit{intervention}. The intervention of vehicle $k$ at station $i$ at time $t$ is denoted as $O_i^k(t)$. The value of $O_i^k(t)$ is reported as follows. $O_i^k(t) > 0$ if bikes were unloaded from the vehicle to the station (i.e., the number of bikes at the station increased), and $O_i^k(t) < 0$ if bikes were loaded on the vehicle. In each epoch $t'$, $O_i^k(t)$ is known for $t \leq t'$ for all stations $i$.
The decisions regarding the interventions of the vehicles involves dynamic non-deterministic elements, i.e., the stock of the station at a future arrival time. The decision on the intervention as well as the visit of a station may therefore be revised at reoptimization epochs as the information on the stock of the station is updated, either because the planned intervention is now infeasible, e.g., the load of bikes was planned but the station is empty as the vehicle arrives at the station, or because the planned intervention is now sub-optimal, e.g., visiting another station is now more beneficial to the quality of service.

Decisions on the routing of the vehicles and loading and unloading operations are planned to minimize the expected shortages. We define an expected \textit{bike shortage} for a station as the time interval starting when the forecast of the stock indicates an empty station while the balance indicates that there is demand for bike rentals.
Similarly, we define an expected \textit{stand shortage} as the time interval starting when the forecast of the stock indicates a full station while the balance indicates that there is demand for bike returns.

Figure \ref{fig:trend} (top) reports an example of the quantities described in this section. In the example, a bike shortage is forecasted in the interval $(t_1, t_2)$ and a stand shortage is forecasted in interval $(t_3,t_4)$. In Figure \ref{fig:trend} (bottom) the example is replicated with the addition of one intervention at time $t''$, where $2$ bikes are removed from station $i$ and loaded on vehicle $k$.

\subsubsection{Optimization algorithm}

\begin{figure}[H]
	\centering
	\includegraphics[width=1.0\linewidth]{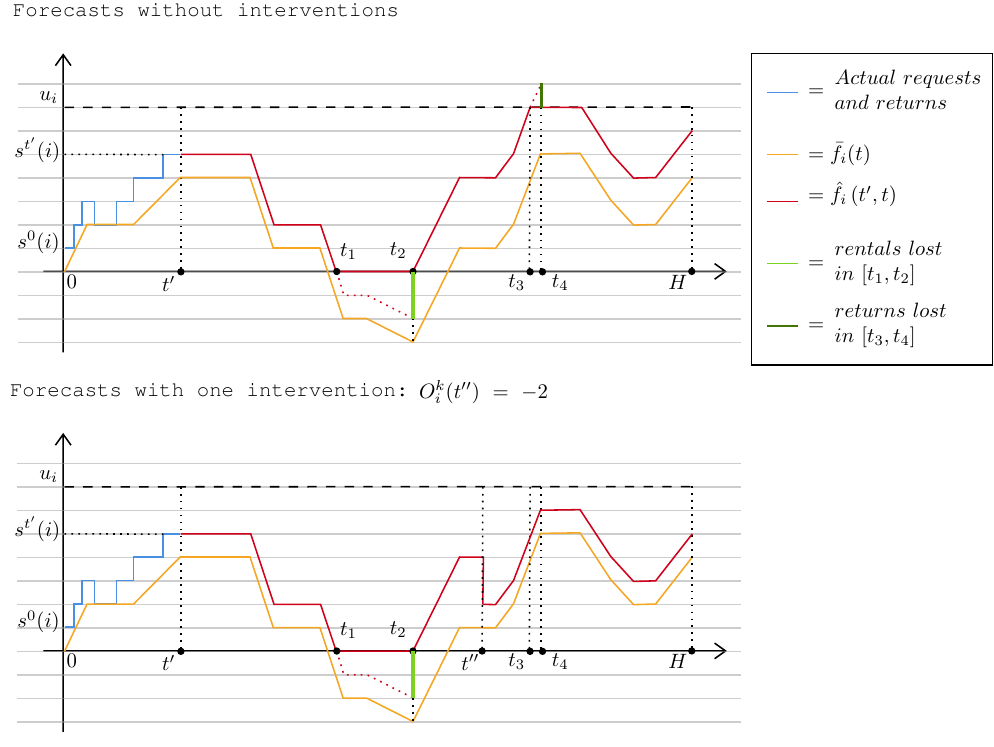}
	\caption{Example of station $i$ during a simulation with zero interventions (top) and with one intervention (bottom) scheduled.}
	\label{fig:trend}
\end{figure}

The proposed algorithm for the solution of the relocation problem in each reoptimization epoch is reported in Algorithm \ref{algo:reopt}. The algorithm is based on the concept of forecasted shortage. As stated above, a shortage is defined as a time interval where the station has either no vacant stands and the balance forecast for the station indicates that there are requests for the return of bikes (stand shortage) or no bike is available for rental and the balance forecast indicates that there would be requests for rental (bike shortage). To forecast shortages it is therefore necessary to forecast the stock of each station. At time $t'$, the forecast of the stock of station $i$ at time $t > t'$, i.e., the forecast of the number of bikes parked at the station at each time $t > t'$, can be computed given the forecast of the balance, $\bar{f}_{i}(t)$, and the stock at time $t'$, $s_i(t')$. The forecast of the stock is denoted as $\hat{f}_{i}(t',t)$ and is computed as follows.
First, $\bar{f}_{i}(t)$ is shifted by $s_i(t') - \bar{f}_{i}(t')$. Then, the resulting function is censored when it exceeds $0$ as a lower bound and $u_i$ as an upper bound (see Figure \ref{fig:trend}, top). A censoring of the function for the lower bound is exemplified in Figure \ref{fig:trend} (top) in the $[t_1, t_2]$ time interval. In the example, the censoring ends at time $t_2$ as, at that time, the station is forecasted to be empty and the forecast of balance indicates that users will be returning bikes and, thus, the forecast of the stock increases. It must be noted that, because of its definition and interpretation, function $\hat{f}_{i}(t',t)$ is only defined for $t > t'$. Furthermore, the function $\hat{f}_{i}(t',t)$ is also influenced by the interventions that are scheduled to take place in station $i$. An example of the effect of an intervention on the forecast of the stock of a station is shown in Figure \ref{fig:trend} (bottom).

The definition of shortage allows us to define the expected number of lost rentals or returns during the shortage, i.e., the number of rentals that are lost during a bike shortage and the number of returns that are lost during a stand shortage, respectively. Intuitively, shortages with starting time close to the reoptimization time and a larger amount of unsatisfied demand are more urgent than those with starting time further away in time and fewer lost requests. Because of the dynamic and stochastic nature of the information, shortages that are further away in time are also less certain in their definition than those that are closer to the current time. This allows us to define a measure of priority for serving stations.

At any time, all vehicles whose shift has not ended (i.e., vehicles that are on duty or whose shift for the day is yet to start), are either at the depot or performing bike relocations. In the former case, the vehicle is loaded with half its capacity of bikes. In the latter case, the vehicle is either at a station or en route to a station. In both cases, at reoptimization epochs, the operations to be performed at such a station are reassessed based on up-to-date information on the status of the system (see Algorithm \ref{algo:reopt}, Line \ref{algo:r2}) so that the next shortage is postponed as much as possible.

In Line \ref{algo:r3}, the remaining stations are sorted according to a criterion considering the urgency to visit each of the stations, defined based on the time of first shortage and the amount of lost requests. The sorting criterion takes a parameter measuring the alertness to the risk of shortage of the station. In particular, a parameter $\Delta$ is defined such that a station $i$ is considered to be suffering from a bike shortage if the level is equal to $\Delta$ and the balance indicates that there is demand for bike rentals, and suffering from a stand shortage if the level is equal to $u_i - \Delta$ and there is demand for bike returns. Different criteria and values for the parameters for each criterion have been tested. The results are reported in Section \ref{sec:tuning}.

The algorithm goes over the list of stations ordered by non-increasing priority and tries to schedule a visit to the station under exam to the route of the vehicle that is expected to reduce the number of lost requests the most (see Lines \ref{algo:r5} and \ref{algo:r6}).
The ideal intervention is the one that postpones the shortage of the station the most, ideally eliminating any shortage from the station for the day. The arrival time of a vehicle to a station and the number of bikes that are loaded/unloaded are tightly connected. Arriving too early means that there might be not enough bikes to be loaded or vacant stands to unload bikes to. Arriving too late might mean that some rentals/returns have already been lost and the intervention is not up to date with the status and the forecast of the station. In this case, performing the intervention as scheduled might do more harm than good (e.g., a late delivery of bikes may cause a stand shortage at a later time). The intervention is also connected to the vehicle status. A vehicle might have the necessary amount of bikes/free space for the required intervention but be too far to arrive where needed at a useful time. Conversely, a vehicle may be readily available but lack the required number of bikes/free space. Compromise solutions, e.g., suboptimal intervention, are therefore considered. For instance, a vehicle might be scheduled to visit a station to mitigate a bike shortage even if its bike load is not enough to solve the shortage. In this case, however, the best arrival time of the vehicle at the station might differ from the one previously computed, as fewer bikes will be unloaded.

\begin{algorithm}[H]
	\small
	\caption{} \label{algo:reopt}
	\begin{algorithmic}[1]
		\STATE \textbf{Input:} Current time. Vehicles: current load and location. Stations: current stock and balance forecast
		\STATE Re-evaluate actions of each vehicle at current or destination stations \label{algo:r2}
		\STATE Sort remaining stations according to the selected criterion and alertness parameter \label{algo:r3}
		\FOR {\textbf{each} station}
		\STATE Find best vehicle to reduce total amount of shortage \label{algo:r5}
		\STATE Add the station to the tour of the vehicle \label{algo:r6}
		\ENDFOR
	\end{algorithmic}
\end{algorithm}

It must be noted that reoptimizations are performed without allowing detours, meaning that, if a vehicle is en route to a station, only its tour after the visit of such station can be modified.

\section{A case study}\label{sec:casestudy}

\begin{figure}[h]
	\centering
	\includegraphics[width=0.7\linewidth]{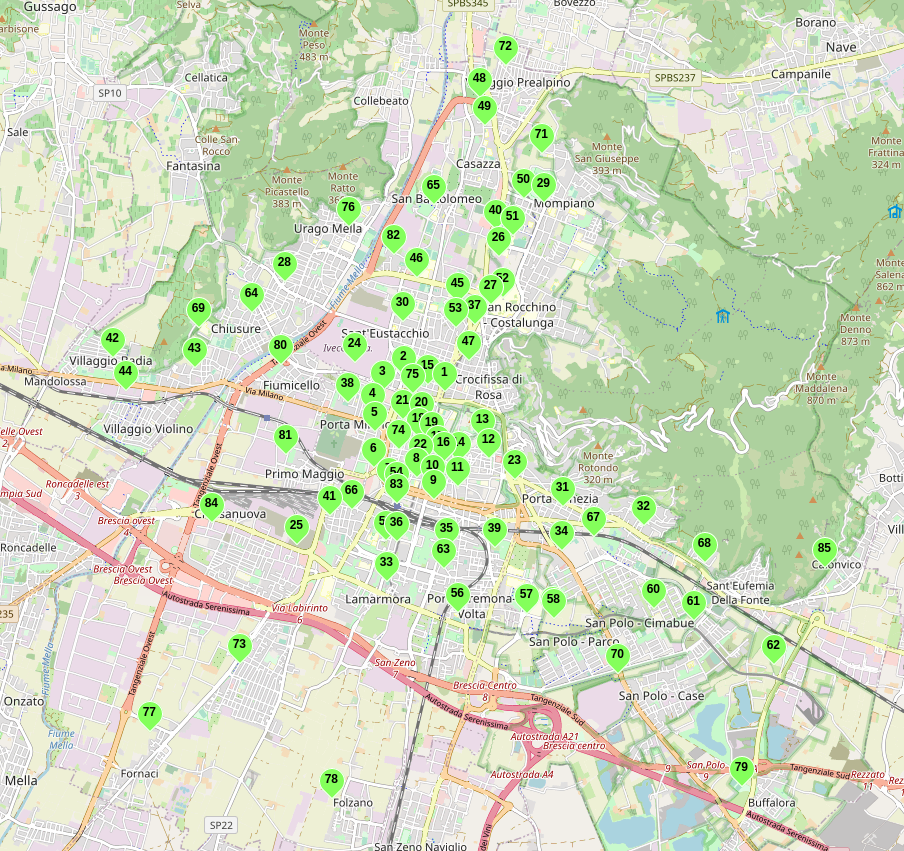}
	\caption{The geographical layout of the 86 station of the Bicimia service.}
	\label{fig:map}
\end{figure}

\begin{figure}[h]
	\centering
	\includegraphics[width=0.7\linewidth]{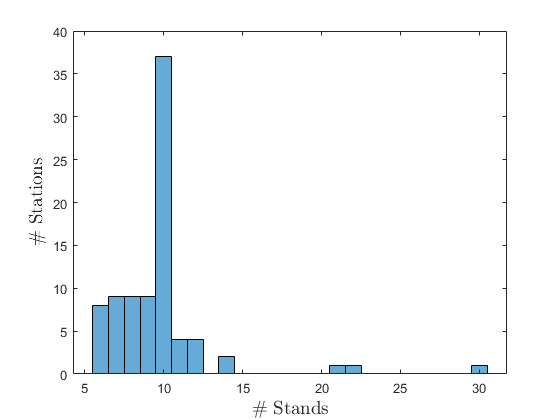}
	\caption{Histogram of the distribution of the capacity of the stations.}
	\label{fig:histo_size}
\end{figure}

The simulation framework presented in this paper has been developed in collaboration with Brescia Mobilità. The company runs a station-based BSS called Bicimia, serving the municipality of Brescia, Italy, and a few neighboring towns. In its current status, the service is characterized by 86 stations (85 bike stations and one depot) covering an area of about 55 $km^2$ (13600 acres). See Figure \ref{fig:map} for the geographical layout of the stations. The stations have an average capacity of about 10 stands with a minimum of 6 and a maximum of 30. The modal size is 10 stands, with 37 stations having such capacity. A total of 7 stations have 6 stands. These stations are generally located in the suburbs of the city. The only station with 30 stands is located just outside the train station. On average, 827 stands are installed in the system. Typically, 360 bikes are present in the system, either parked at a stand or rented. A histogram of the distribution of the capacity of the stations is shown in Figure \ref{fig:histo_size}.

There are currently about 30000 users subscribed to the service. The number of daily trips fluctuates between about 500, e.g., on a rainy winter Sunday, to up to 3000, e.g., on a sunny working day of May. An average of 450 bikes are available for rental and the average rental period is 24 minutes. As opposed to other BSSs (see, for instance, \cite{dell2014bike}), the service is available 24/7. The shifts of the operators cover the time interval from 7:00 to 19:30 in the working days (with a morning shift from 7:00 to 15:00 and an afternoon shift from 11:30 to 19:30) and 7:00 to 13:00 in the pre-festive days (i.e., on Saturdays). Generally, festive days are such also for the operators. Currently, Brescia Mobilità has a heterogeneous fleet of three vehicles. Two vehicles are vans, holding up to 14 bikes. The third vehicle is a flatbed van with a rack holding up to 20 bikes. Typically, in working days, two vehicles are scheduled in the morning shift and one in the afternoon, while in pre-festive days one vehicle is scheduled for the shift. The workforce of Brescia Mobilità for Bicimia and the fleet of vehicles goes beyond the three operators and vehicles generally scheduled for the relocation of bikes. In fact, it is not rare that additional operators, usually performing different tasks (e.g., bike repair), are deployed to perform bike relocation. Also, additional vehicles may be used which are normally assigned to other services. For this reason, in the following, different settings will be tested without explicitly considering the economical implications of employing additional workforce and vehicles.

Currently, the bike relocation problem is not addressed by the company's decision support system. The operators rely on their experience to decide how to relocate the bikes during the day.
As no codified knowledge is available for the forecast of rentals and returns for the Bicimia service, the characteristics of a day influencing the forecast of the requests have been identified according to the experience of the operators. Such characteristics are the following:

\begin{itemize}
	\item type of day, i.e., working day, Saturday, or Sunday;
	\item month of the year;
	\item weather, i.e., sunny or rainy.
\end{itemize}

These characteristics appear to be consistent with those reported in \cite{eren2020review}, where factors affecting bike-sharing demand are reviewed.
For each combination of these characteristics, a past day has been identified to represent the forecast of bike rentals and returns in each station and, thus, functions $\bar{f}_{i}(t)$. While using data of successful pickups and deliveries of bikes might be the cause of bias in functions $\bar{f}_{i}(t)$,
the use of the simulator framework, together with other DSS tools currently being developed for the company, will mitigate (and hopefully eliminate) this effect by reducing customer dissatisfaction.

\section{Computational results} \label{sec:compres}

In this section we report how instances to test the simulation framework have been generated.
A set of experiments is then presented to showcase the decisions that the proposed framework can support.
As the framework is aimed at improving the quality of the existing service, the baseline for every experiment is the current layout and fleet of the service as described in Section \ref{sec:casestudy}.
For the sake of ease of interpretation, all vehicles are assigned a capacity equal to 14 bikes. This also has a practical justification: the company wishes to make its fleet more homogeneous to reduce maintenance costs.
Results are reported for a sunny working day of May, that is, the type of day with the highest usage of the system. In all cases, results are presented as the average over 100 scenarios, with an average run time of around 4 seconds per scenario.

The stochastic processes from which we generated the user request scenarios have been introduced in Section \ref{sec:framework}.
The parameters $\lambda_{i}^+(t)$ and $\lambda_{i}^-(t)$ that characterize the processes at station $i$ are determined according to data provided by the company.
In particular, the company identified, for each type of day, a single date which is representative of the average use of the system.
The stepwise cumulative curves of bike returns (rentals) detected by the BSS in each station were then approximated by a simple Piecewise Linear function (PL).
The approximation procedure initializes the PL function with one segment representing the overall trend of the day.
The time of the day where the PL and the cumulative curve are the most distant is then evaluated.
If the distance is greater than a specified threshold a new breakpoint for the PL is created, thus splitting the corresponding interval in two.
This step is iterated until PL is reasonably close to the cumulative curve.
In general, a small number of iterations are enough to meet the stopping condition and obtain a simple PL whose slope define $\lambda_{i}^+(t)$ ($\lambda_{i}^-(t)$) for the day.
Scenarios are then sampled from the resulting stochastic processes.
The forecast of the balance function $\bar{f}_i(t)$ is obtained in a similar way starting from the difference of the cumulative curves of rentals and returns.

\medskip

The framework has been implemented in Java running on a Windows 10 machine equipped with a 3.5 GHz Intel Xeon E5-1650v2 processor and 64GB of RAM.

The remainder of this section is structured as follows. In Section \ref{sec:tuning} a discussion on the tuning of two parameters of the DBRP algorithm is presented. In Section \ref{sec:cs}, results on the system in its current status are presented. Then, experiments on the shifts of the vehicles are reported in Section \ref{sec:shifts}, experiments on the composition of the fleet of vehicles are reported in Section \ref{sec:fleetv}, experiments on the sizing of the stations are reported in Section \ref{sec:sizes}, and experiments on the sizing of the fleet of bikes are reported in Section \ref{sec:bikefleet}.

\subsection{Parameter selection} \label{sec:tuning}

Different criteria for the sorting of unvisited stations and values for the alertness parameter $\Delta$ have been tested. With respect to the sorting criteria, the following have been tested:
\begin{itemize}
	\item First Shortage Time (FST): stations are sorted based on the starting time of the next shortage.
	\item Group by Shortage Time + Amount (GSTA): First, stations are sorted based on the starting time of the first shortage. Then, stations are grouped together if such time happens in the same quarter of an hour (e.g., stations with shortages happening between 12:00 and 12:15 are grouped together). Within a group, stations with the higher number of expected lost rentals or returns have higher priority.
	\item Group by Shortage Time + Rate (GSTR): The sorting is made similarly to the GSTA criterion but, within a group, stations are sorted based on the number of expected lost rentals or returns divided by the expected duration of the shortage (i.e. higher priority is given to station loosing bikes faster).
\end{itemize}

\begin{figure}[H]
	\centering
	\includegraphics[width=0.8\linewidth]{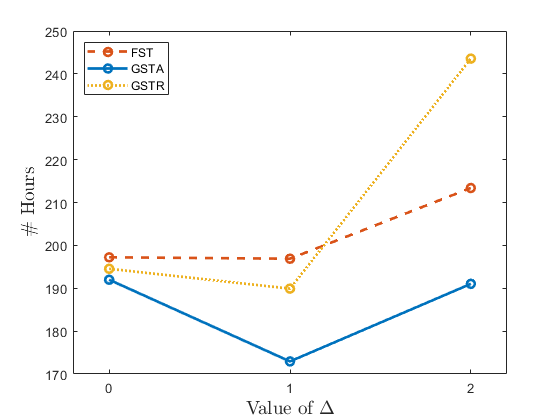}
	\caption{The performance of the system for different sorting criteria and values of $\Delta$.}
	\label{fig:tuning}
\end{figure}

The results of the selection of sorting criterion and alertness parameter are reported in Figure \ref{fig:tuning}, where the x-axis reports the values of $\Delta \in \{0,1,2\}$ and the y-axis the performance of the system. For each criterion-parameter combination, the performance is measured as the number of hours where stations are either full or empty. The result is reported as the average across different combinations for the number of vehicles allocated in the morning and afternoon shift, up to a maximum of six vehicles across the two shifts (a more in depth discussion on these results is presented in Section \ref{sec:numopt} for the best performing combination). Results indicate that the best performance is obtained by the GSTA criterion with $\Delta = 1$. Results presented hereafter are based on this combination.

\subsection{Performance of the system with no relocation} \label{sec:cs}

\begin{figure}[H]
	\centering
	\includegraphics[width=0.32\linewidth]{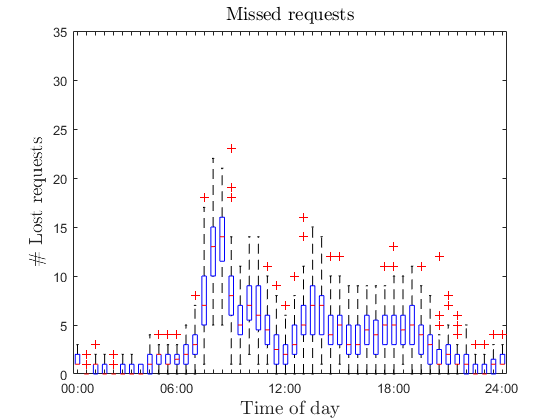}
	\includegraphics[width=0.32\linewidth]{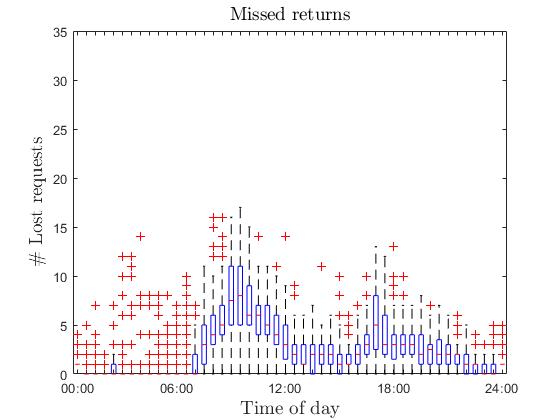}
	\includegraphics[width=0.32\linewidth]{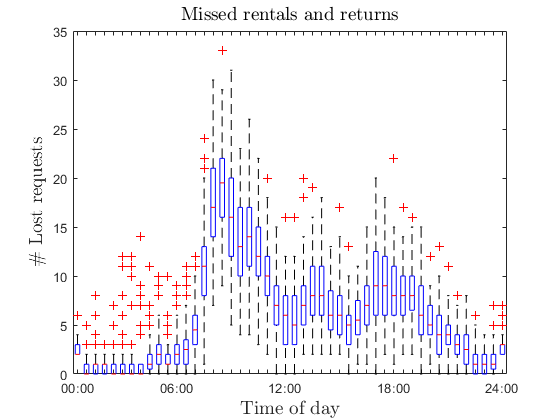}
	\caption{The number of missed rentals (left), returns (center), and the sum of the two (right) over the generated scenarios when no vehicle is performing the relocation.}
	\label{fig:missed_noveicles}
\end{figure}

To show the benefits of an accurate sizing of the fleet and scheduling of the operators, the performance of the system when no vehicle is in service is reported as simulated by the presented framework. In Figure \ref{fig:missed_noveicles} the distribution of the number of missed rentals, returns, and the sum of the two across the 100 scenarios is reported in intervals of 30 minutes.
In Figure \ref{fig:missed_noveicles} (right) one major peak in the number of missed rentals and returns is observed in the morning, in the $[7:30 - 11:30]$ interval, and two minor peaks are observed around the end of high-school classes, i.e., $[13:00 - 14:30]$, and at the end of the typical working day, i.e., $[17:00 - 19:30]$.

\subsection{Experiments on shifts} \label{sec:shifts}

\subsubsection{Number of vehicles and operators} \label{sec:numopt}

In this section, an increasing number of vehicles for each of the shifts that are currently set by the company for working days (i.e., Monday to Friday) is tested. We measure the performance of each tested setting with the average number of hours that stations have been empty or full. Tables \ref{tab:morningEvening_empty}, \ref{tab:morningEvening_full}, and \ref{tab:morningEvening_tot} report the number of hours stations are empty, full, and the sum of the two for each combination of the number of vehicles allocated to each shift, up to a maximum of 6 vehicles, with each vehicle having a capacity of 14 bikes.

\begin{table}[H]
	\centering
	\begin{tabular}{c|r|rrrrrrr}
		\multicolumn{1}{r}{} & \multicolumn{1}{r}{} & \multicolumn{7}{c}{\# Vehicles - Afternoon shift} \\
		\cmidrule{3-9}    \multicolumn{1}{r}{} & \multicolumn{1}{r}{} & 0     & 1     & 2     & 3     & 4     & 5     & 6 \\
		\cmidrule{3-9}    \multirow{7}[1]{*}{\begin{tabular}[l]{@{}c@{}} \# Vehicles - \\ Morning shift \end{tabular}} & 0     & \cellcolor[rgb]{ .973,  .412,  .42}192.45 & \cellcolor[rgb]{ .984,  .58,  .455}168.45 & \cellcolor[rgb]{ .988,  .682,  .475}153.73 & \cellcolor[rgb]{ .992,  .765,  .49}142.20 & \cellcolor[rgb]{ .996,  .812,  .498}135.59 & \cellcolor[rgb]{ .996,  .839,  .502}131.50 & \cellcolor[rgb]{ .996,  .851,  .506}129.97 \\
				  & 1     & \cellcolor[rgb]{ .984,  .616,  .459}163.43 & \cellcolor[rgb]{ .992,  .757,  .486}143.30 & \cellcolor[rgb]{ 1,  .851,  .506}129.83 & \cellcolor[rgb]{ 1,  .898,  .514}123.02 & \cellcolor[rgb]{ .996,  .918,  .514}119.38 & \cellcolor[rgb]{ .957,  .906,  .514}117.59 &  \\
				  & 2     & \cellcolor[rgb]{ .992,  .776,  .49}140.73 & \cellcolor[rgb]{ 1,  .898,  .514}122.91 & \cellcolor[rgb]{ .906,  .894,  .51}115.47 & \cellcolor[rgb]{ .776,  .855,  .502}109.74 & \cellcolor[rgb]{ .698,  .835,  .498}106.51 &       &  \\
				  & 3     & \cellcolor[rgb]{ 1,  .922,  .518}119.42 & \cellcolor[rgb]{ .718,  .839,  .498}107.28 & \cellcolor[rgb]{ .576,  .8,  .49}101.21 & \cellcolor[rgb]{ .514,  .78,  .486}98.44 &       &       &  \\
				  & 4     & \cellcolor[rgb]{ .776,  .855,  .502}109.86 & \cellcolor[rgb]{ .506,  .776,  .486}98.07 & \cellcolor[rgb]{ .404,  .749,  .482}93.66 &       &       &       &  \\
				  & 5     & \cellcolor[rgb]{ .643,  .816,  .494}104.05 & \cellcolor[rgb]{ .388,  .745,  .482}92.96 &       &       &       &       &  \\
				  & 6     & \cellcolor[rgb]{ .569,  .796,  .49}100.81 &       &       &       &       &       &  \\
			\end{tabular}
		
	\caption{The average number of hours that stations are empty, for each combination of morning and afternoon shifts.}
	\label{tab:morningEvening_empty}
 \end{table}

\begin{table}[H]
	\centering
	\begin{tabular}{c|r|rrrrrrr}
		\multicolumn{1}{r}{} & \multicolumn{1}{r}{} & \multicolumn{7}{c}{\# Vehicles - Afternoon shift} \\
		\cmidrule{3-9}    \multicolumn{1}{r}{} & \multicolumn{1}{r}{} & 0     & 1     & 2     & 3     & 4     & 5     & 6 \\
		\cmidrule{3-9}    \multirow{7}[1]{*}{\begin{tabular}[l]{@{}c@{}} \# Vehicles - \\ Morning shift \end{tabular}} & 0     & \cellcolor[rgb]{ .973,  .412,  .42}81.38 & \cellcolor[rgb]{ .988,  .671,  .471}67.32 & \cellcolor[rgb]{ .992,  .749,  .486}62.93 & \cellcolor[rgb]{ .996,  .824,  .502}58.75 & \cellcolor[rgb]{ .996,  .851,  .506}57.32 & \cellcolor[rgb]{ 1,  .898,  .514}54.76 & \cellcolor[rgb]{ .992,  .918,  .514}53.22 \\
				  & 1     & \cellcolor[rgb]{ .984,  .612,  .459}70.43 & \cellcolor[rgb]{ .996,  .82,  .498}58.93 & \cellcolor[rgb]{ 1,  .922,  .518}53.38 & \cellcolor[rgb]{ .851,  .878,  .506}50.72 & \cellcolor[rgb]{ .745,  .847,  .502}48.90 & \cellcolor[rgb]{ .608,  .808,  .494}46.46 &  \\
				  & 2     & \cellcolor[rgb]{ .992,  .722,  .482}64.34 & \cellcolor[rgb]{ .988,  .918,  .514}53.12 & \cellcolor[rgb]{ .694,  .831,  .498}47.94 & \cellcolor[rgb]{ .553,  .792,  .49}45.53 & \cellcolor[rgb]{ .388,  .745,  .482}42.58 &       &  \\
				  & 3     & \cellcolor[rgb]{ .992,  .749,  .486}62.87 & \cellcolor[rgb]{ .886,  .886,  .51}51.33 & \cellcolor[rgb]{ .608,  .808,  .494}46.45 & \cellcolor[rgb]{ .4,  .749,  .482}42.84 &       &       &  \\
				  & 4     & \cellcolor[rgb]{ .992,  .765,  .49}62.05 & \cellcolor[rgb]{ .824,  .871,  .506}50.22 & \cellcolor[rgb]{ .502,  .776,  .486}44.62 &       &       &       &  \\
				  & 5     & \cellcolor[rgb]{ .996,  .78,  .494}61.09 & \cellcolor[rgb]{ .749,  .847,  .502}48.95 &       &       &       &       &  \\
				  & 6     & \cellcolor[rgb]{ .996,  .8,  .494}60.07 &       &       &       &       &       &  \\
			\end{tabular}
	\caption{The average number of hours that stations are full for each combination of morning and afternoon shifts.}
	\label{tab:morningEvening_full}
 \end{table}

\begin{table}[H]
 	\centering
	 \begin{tabular}{c|r|rrrrrrr}
		\multicolumn{1}{r}{} & \multicolumn{1}{r}{} & \multicolumn{7}{c}{\# Vehicles - Afternoon shift} \\
		\cmidrule{3-9}    \multicolumn{1}{r}{} & \multicolumn{1}{r}{} & 0     & 1     & 2     & 3     & 4     & 5     & 6 \\
		\cmidrule{3-9}    \multirow{7}[1]{*}{\begin{tabular}[l]{@{}c@{}} \# Vehicles - \\ Morning shift \end{tabular}} & 0     & \cellcolor[rgb]{ .973,  .412,  .42}273.83 & \cellcolor[rgb]{ .984,  .604,  .459}235.77 & \cellcolor[rgb]{ .988,  .702,  .478}216.66 & \cellcolor[rgb]{ .996,  .78,  .494}200.95 & \cellcolor[rgb]{ .996,  .824,  .502}192.91 & \cellcolor[rgb]{ 1,  .855,  .506}186.26 & \cellcolor[rgb]{ 1,  .871,  .51}183.19 \\
				  & 1     & \cellcolor[rgb]{ .984,  .616,  .459}233.86 & \cellcolor[rgb]{ .992,  .776,  .49}202.23 & \cellcolor[rgb]{ 1,  .871,  .51}183.21 & \cellcolor[rgb]{ 1,  .918,  .518}173.74 & \cellcolor[rgb]{ .918,  .898,  .51}168.28 & \cellcolor[rgb]{ .843,  .875,  .506}164.05 &  \\
				  & 2     & \cellcolor[rgb]{ .992,  .761,  .49}205.07 & \cellcolor[rgb]{ 1,  .906,  .518}176.03 & \cellcolor[rgb]{ .831,  .871,  .506}163.41 & \cellcolor[rgb]{ .686,  .831,  .498}155.27 & \cellcolor[rgb]{ .576,  .8,  .49}149.09 &       &  \\
				  & 3     & \cellcolor[rgb]{ 1,  .875,  .51}182.29 & \cellcolor[rgb]{ .745,  .847,  .502}158.61 & \cellcolor[rgb]{ .553,  .792,  .49}147.66 & \cellcolor[rgb]{ .439,  .757,  .482}141.28 &       &       &  \\
				  & 4     & \cellcolor[rgb]{ .98,  .914,  .514}171.91 & \cellcolor[rgb]{ .565,  .796,  .49}148.29 & \cellcolor[rgb]{ .388,  .745,  .482}138.28 &       &       &       &  \\
				  & 5     & \cellcolor[rgb]{ .863,  .878,  .506}165.14 & \cellcolor[rgb]{ .451,  .761,  .482}141.91 &       &       &       &       &  \\
				  & 6     & \cellcolor[rgb]{ .788,  .859,  .502}160.88 &       &       &       &       &       &  \\
			\end{tabular}
 	\caption{The average number of hours that stations are either full or empty, for each combination of morning and afternoon shifts.}
	\label{tab:morningEvening_tot}
\end{table}

\begin{figure}[H]
	\centering
	\includegraphics[width=0.3\linewidth]{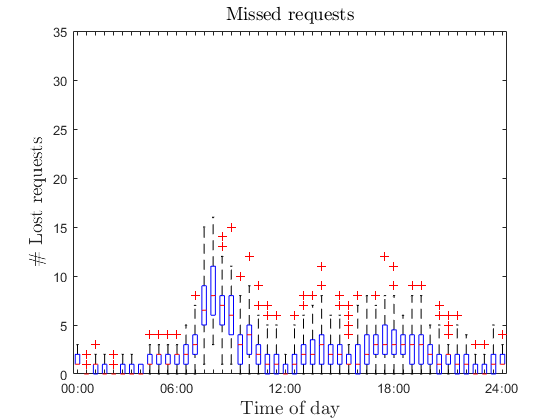}
	\includegraphics[width=0.3\linewidth]{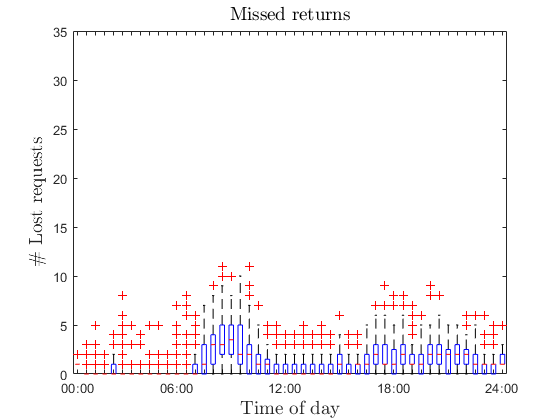}
	\includegraphics[width=0.3\linewidth]{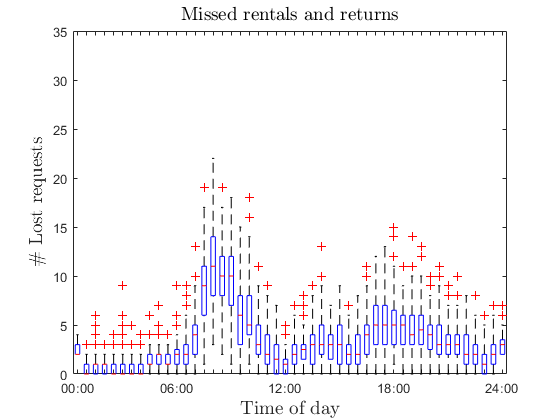}
	\caption{The number of missed rentals, returns, and the sum of the two over the generated scenarios with two vehicles in the morning shift and one in the afternoon shift.}
	\label{fig:missed_2m1a}
\end{figure}

Several insights can be gained from Table \ref{tab:morningEvening_tot}. For instance, given that the company currently employs three operators, the simulator finds that the best allocation of operators and vehicles is to have two vehicles in the morning shift and one in the afternoon shift. This validates the current practice of the company. It should be noted, however, that the company does not currently reposition bikes according to the algorithm. The effectiveness of this scheduling can be measured by analyzing Figure \ref{fig:missed_2m1a} in comparison with Figure \ref{fig:missed_noveicles}. The benefit of bike relocation on the number of lost rentals and returns is clear, especially during and after the [11:30 - 15:00] interval, when all three vehicles are active. The number of lost rentals and returns in the morning peak also has decreased, albeit to a lesser degree, in particular in the initial (and more intense) part of the peak. Concerning the number of hours that any station is either full or empty, the current schedule is simulated  to result in 176.03 hours, a reduction of about 36\% over the 273.83 hours from the simulation when the relocation of bikes is not performed.

Given a certain allocation to shifts, the simulator allows for the assessment of the impact of the addition or removal of one operator/vehicle. For instance, if an  operator/vehicle was to be added to the current roster, it would be best to schedule it in the morning shift, with a reduction of the total number of hours that stations are either full or empty of about 9.8\%.

By comparing the results in Tables \ref{tab:morningEvening_empty} and \ref{tab:morningEvening_full}, it can be observed that the hours that stations are empty greatly exceed the hours that stations are full. This could indicate that the fleet of bikes in the system is possibly undersized. Experiments regarding the sizing of the fleet of bikes are reported in Section \ref{sec:bikefleet}.

\subsubsection{No overlaps} \label{sec:nooverlap}

In the current status of the system, three vehicles are required to cover the shifts of the operators. This is because the morning and afternoon shifts are overlapping and therefore three operators are on schedule at the same time. The overlap, between 11:30 and 15:00, existing in the current organization of shifts, is justified by the peak of demand that happens around lunchtime because of students exiting schools. The company could consider redesigning the shifts so that no overlap is present. This would allow them to reduce the number of operators working at the same time, and therefore the number of vehicles required to perform the relocation.
The impact of this redesign is assessed by testing the following shifts:

\begin{itemize}
	\item morning shift: from 5:00 to 13:00;
	\item afternoon shift: from 13:00 to 21:00.
\end{itemize}
This redesign also has the advantage of better covering the morning and evening peaks.

\begin{figure}[H]
	  \centering
	  \includegraphics[width=0.7\linewidth]{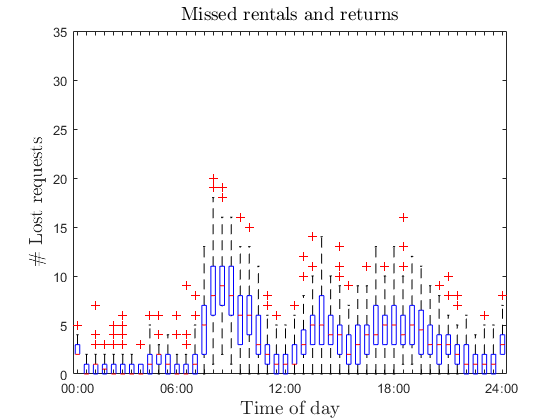}
	  \caption{The number of missed rentals and returns over the generated scenarios with two vehicles in the morning shift and one in the afternoon shift (non-overlapping shifts).}
	  \label{fig:miss_nooverlap}
  \end{figure}

Results are reported in Tables \ref{tab:nool_empty}, \ref{tab:nool_full}, and \ref{tab:nool_total} with the same interpretation of Tables \ref{tab:morningEvening_empty}, \ref{tab:morningEvening_full}, and \ref{tab:morningEvening_tot}, respectively. The different working shifts presented in this section allow for a greater reduction of the number of hours in which stations are either empty or full. This is due to a better timing of the shifts with respect to the peaks in the usage of the system. In Figure \ref{fig:miss_nooverlap} the number of missed rentals and returns is reported for the case of two vehicles scheduled in the morning shift and one scheduled in the afternoon, i.e., the schedule closest to the current status of the system. In comparison with Figure \ref{fig:missed_2m1a} (right), the morning peak is reduced, whereas the noon peak suffers from a slight increase, due to the reduced number of vehicles active in that time frame.

A similar configuration of the shifts has been tested by the company. The empirical results indeed confirmed the reduction of the morning peak at the expense of the noon peak. The increase in the noon peak is the reason this configuration of the shifts has not been adopted by the company.

\begin{table}[H]
	\centering
    \begin{tabular}{c|r|rrrrrrr}
		\multicolumn{1}{r}{} & \multicolumn{1}{r}{} & \multicolumn{7}{c}{\# Vehicles - Afternoon shift} \\
	\cmidrule{3-9}    \multicolumn{1}{r}{} & \multicolumn{1}{r}{} & 0     & 1     & 2     & 3     & 4     & 5     & 6 \\
	\cmidrule{3-9}    \multirow{7}[1]{*}{\begin{tabular}[l]{@{}c@{}} \# Vehicles - \\ Morning shift \end{tabular}} & 0     & \cellcolor[rgb]{ .973,  .412,  .42}192.45 & \cellcolor[rgb]{ .984,  .565,  .451}170.26 & \cellcolor[rgb]{ .988,  .635,  .463}159.62 & \cellcolor[rgb]{ .992,  .71,  .478}148.64 & \cellcolor[rgb]{ .992,  .753,  .486}142.49 & \cellcolor[rgb]{ .996,  .78,  .49}138.45 & \cellcolor[rgb]{ .996,  .788,  .494}137.26 \\
			  & 1     & \cellcolor[rgb]{ .984,  .596,  .455}165.74 & \cellcolor[rgb]{ .992,  .741,  .486}144.13 & \cellcolor[rgb]{ .996,  .812,  .498}133.64 & \cellcolor[rgb]{ 1,  .882,  .51}123.46 & \cellcolor[rgb]{ 1,  .918,  .518}118.46 & \cellcolor[rgb]{ .984,  .918,  .514}116.42 &  \\
			  & 2     & \cellcolor[rgb]{ .992,  .769,  .49}140.25 & \cellcolor[rgb]{ 1,  .898,  .514}121.11 & \cellcolor[rgb]{ .918,  .898,  .51}111.08 & \cellcolor[rgb]{ .824,  .871,  .506}103.55 & \cellcolor[rgb]{ .769,  .855,  .502}99.16 &       &  \\
			  & 3     & \cellcolor[rgb]{ .922,  .898,  .51}111.48 & \cellcolor[rgb]{ .706,  .835,  .498}94.36 & \cellcolor[rgb]{ .608,  .808,  .494}86.52 & \cellcolor[rgb]{ .537,  .788,  .49}80.88 &       &       &  \\
			  & 4     & \cellcolor[rgb]{ .69,  .831,  .498}92.90 & \cellcolor[rgb]{ .482,  .773,  .486}76.47 & \cellcolor[rgb]{ .388,  .745,  .482}68.84 &       &       &       &  \\
			  & 5     & \cellcolor[rgb]{ .592,  .804,  .494}85.34 & \cellcolor[rgb]{ .396,  .745,  .482}69.60 &       &       &       &       &  \\
			  & 6     & \cellcolor[rgb]{ .545,  .788,  .49}81.34 &       &       &       &       &       &  \\
		\end{tabular}
	
	  \caption{The average number of hours that stations are empty, for each combination of morning and afternoon non-overlapping shifts.}
	\label{tab:nool_empty}
  \end{table}

\begin{table}[H]
	\centering
    \begin{tabular}{c|r|rrrrrrr}
		\multicolumn{1}{r}{} & \multicolumn{1}{r}{} & \multicolumn{7}{c}{\# Vehicles - Afternoon shift} \\
	\cmidrule{3-9}    \multicolumn{1}{r}{} & \multicolumn{1}{r}{} & 0     & 1     & 2     & 3     & 4     & 5     & 6 \\
	\cmidrule{3-9}    \multirow{7}[1]{*}{\begin{tabular}[l]{@{}c@{}} \# Vehicles - \\ Morning shift \end{tabular}} & 0     & \cellcolor[rgb]{ .973,  .412,  .42}81.38 & \cellcolor[rgb]{ .988,  .667,  .471}62.25 & \cellcolor[rgb]{ .996,  .784,  .494}53.35 & \cellcolor[rgb]{ 1,  .855,  .506}47.98 & \cellcolor[rgb]{ 1,  .871,  .51}46.76 & \cellcolor[rgb]{ 1,  .878,  .51}46.10 & \cellcolor[rgb]{ 1,  .886,  .514}45.61 \\
			  & 1     & \cellcolor[rgb]{ .984,  .604,  .459}66.89 & \cellcolor[rgb]{ 1,  .875,  .51}46.32 & \cellcolor[rgb]{ .784,  .859,  .502}36.67 & \cellcolor[rgb]{ .643,  .816,  .494}32.70 & \cellcolor[rgb]{ .596,  .804,  .494}31.33 & \cellcolor[rgb]{ .561,  .792,  .49}30.37 &  \\
			  & 2     & \cellcolor[rgb]{ .988,  .675,  .471}61.55 & \cellcolor[rgb]{ .894,  .89,  .51}39.73 & \cellcolor[rgb]{ .557,  .792,  .49}30.24 & \cellcolor[rgb]{ .459,  .765,  .486}27.41 & \cellcolor[rgb]{ .42,  .753,  .482}26.32 &       &  \\
			  & 3     & \cellcolor[rgb]{ .988,  .667,  .471}62.07 & \cellcolor[rgb]{ .835,  .871,  .506}38.09 & \cellcolor[rgb]{ .525,  .784,  .49}29.34 & \cellcolor[rgb]{ .388,  .745,  .482}25.40 &       &       &  \\
			  & 4     & \cellcolor[rgb]{ .988,  .663,  .471}62.54 & \cellcolor[rgb]{ .847,  .875,  .506}38.42 & \cellcolor[rgb]{ .486,  .773,  .486}28.28 &       &       &       &  \\
			  & 5     & \cellcolor[rgb]{ .988,  .659,  .471}62.82 & \cellcolor[rgb]{ .894,  .89,  .51}39.83 &       &       &       &       &  \\
			  & 6     & \cellcolor[rgb]{ .988,  .671,  .471}61.79 &       &       &       &       &       &  \\
		\end{tabular}
	
	  \caption{The average number of hours that stations are full, for each combination of morning and afternoon non-overlapping shifts.}
	  \label{tab:nool_full}
  \end{table}

\begin{table}[H]
	\centering
    \begin{tabular}{c|r|rrrrrrr}
		\multicolumn{1}{r}{} & \multicolumn{1}{r}{} & \multicolumn{7}{c}{\# Vehicles - Afternoon shift} \\
	\cmidrule{3-9}    \multicolumn{1}{r}{} & \multicolumn{1}{r}{} & 0     & 1     & 2     & 3     & 4     & 5     & 6 \\
	\cmidrule{3-9}    \multirow{7}[1]{*}{\begin{tabular}[l]{@{}c@{}} \# Vehicles - \\ Morning shift \end{tabular}} & 0     & \cellcolor[rgb]{ .973,  .412,  .42}273.83 & \cellcolor[rgb]{ .984,  .592,  .455}232.51 & \cellcolor[rgb]{ .988,  .678,  .471}212.97 & \cellcolor[rgb]{ .992,  .749,  .486}196.62 & \cellcolor[rgb]{ .996,  .78,  .49}189.25 & \cellcolor[rgb]{ .996,  .8,  .494}184.55 & \cellcolor[rgb]{ .996,  .808,  .498}182.87 \\
			  & 1     & \cellcolor[rgb]{ .984,  .592,  .455}232.63 & \cellcolor[rgb]{ .992,  .773,  .49}190.45 & \cellcolor[rgb]{ 1,  .863,  .506}170.31 & \cellcolor[rgb]{ 1,  .922,  .518}156.16 & \cellcolor[rgb]{ .937,  .902,  .514}149.79 & \cellcolor[rgb]{ .906,  .894,  .51}146.79 &  \\
			  & 2     & \cellcolor[rgb]{ .992,  .725,  .482}201.80 & \cellcolor[rgb]{ 1,  .902,  .514}160.84 & \cellcolor[rgb]{ .847,  .875,  .506}141.32 & \cellcolor[rgb]{ .737,  .843,  .502}130.96 & \cellcolor[rgb]{ .682,  .827,  .498}125.48 &       &  \\
			  & 3     & \cellcolor[rgb]{ .996,  .847,  .506}173.55 & \cellcolor[rgb]{ .753,  .851,  .502}132.45 & \cellcolor[rgb]{ .58,  .8,  .49}115.86 & \cellcolor[rgb]{ .482,  .773,  .486}106.28 &       &       &  \\
			  & 4     & \cellcolor[rgb]{ .996,  .918,  .514}155.44 & \cellcolor[rgb]{ .573,  .796,  .49}114.89 & \cellcolor[rgb]{ .388,  .745,  .482}97.12 &       &       &       &  \\
			  & 5     & \cellcolor[rgb]{ .918,  .898,  .51}148.16 & \cellcolor[rgb]{ .514,  .78,  .486}109.43 &       &       &       &       &  \\
			  & 6     & \cellcolor[rgb]{ .867,  .882,  .51}143.13 &       &       &       &       &       &  \\
		\end{tabular}
	  \caption{The average number of hours that stations are either empty or full, for each combination of morning and afternoon non-overlapping shifts.}
	  \label{tab:nool_total}
  \end{table}

\subsubsection{Short shifts} \label{sec:parttime}

As noted in Section \ref{sec:nooverlap}, the overlap in the current shifts of the vehicles is justified by the need to cover the noon peak. For the day under consideration, however, the highest peak happens in the morning, when students and commuters have to reach schools and working places, respectively, and the second highest in the evening.

To test whether it would be beneficial to treat in a similar way the morning and evening peaks of workers commuting to and from work, we investigate a setting with work shifts of the following two kinds only:

\begin{itemize}
	\item morning shift: from 6:00 to 10:00,
	\item evening part-time shift: from 17:00 to 21:00,
\end{itemize}

and assume that full-time operators are available to cover two shifts of 4 hours each and also, if needed, part-time operators are involved to cover one 4-hour shift.

\begin{figure}[H]
	\centering
	\includegraphics[width=0.45\linewidth]{resources/figures/missed/mTOT_2m1a.png}
	\includegraphics[width=0.45\linewidth]{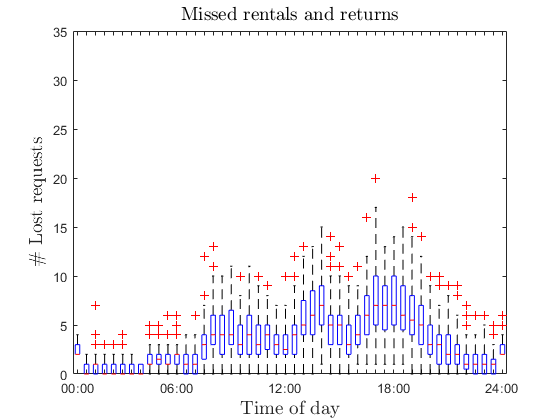}
	\caption{The number of missed rentals and returns over the generated scenarios with full-time shifts (left, two vehicles in the morning, one in the afternoon) and part-time shifts (right, five vehicles in the morning, one in the evening).}
	\label{fig:comparisonparttime}
\end{figure}

Results reported in Table \ref{tab:parttime} provide a valuable insight when compared with those reported in Table \ref{tab:morningEvening_tot}. In the current situation, vehicles are scheduled for a total of 24 hours (i.e., 3 shifts of 8 hours). With the current organization of shifts, the simulator reports that 176.03 hours of either full or empty station are expected. Organizing the work in 4 hours shifts while maintaining a total of 24 man-hour schedule, as presented in this section, could result in 167.42 hours of disservice, a reduction of about 5\% with respect to the current situation, with 5 vehicles scheduled in the morning and 1 in the evening part-time shifts. A comparison between this case and the current schedule of the vehicles is presented in Figure \ref{fig:comparisonparttime}. Should the number of vehicles play a more important role in the decision-making process, since there is no shift overlapping, the scheduling using a maximum of three vehicles would be the one where three vehicles are scheduled for both the morning and evening shifts. In this case, 178.38 hours of stations either empty or full are forecasted, a slight increase over the current status of the system, possibly due to the fact that relocating bikes late in the day has less effect on the system, as demand for bikes fades in the late hours of the day.

\begin{table}[H]
	\centering
	\begin{tabular}{c|r|rrrrrrr}
		\multicolumn{1}{r}{} & \multicolumn{1}{r}{} & \multicolumn{7}{c}{\# Vehicles - Evening shift} \\
	\cmidrule{3-9}    \multicolumn{1}{r}{} & \multicolumn{1}{r}{} & 0     & 1     & 2     & 3     & 4     & 5     & 6 \\
	\cmidrule{3-9}    \multirow{7}[1]{*}{\begin{tabular}[l]{@{}c@{}} \# Vehicles - \\ Morning shift \end{tabular}} & 0     & \cellcolor[rgb]{ .973,  .412,  .42}273.83 & \cellcolor[rgb]{ .98,  .549,  .447}251.32 & \cellcolor[rgb]{ .984,  .62,  .463}240.08 & \cellcolor[rgb]{ .988,  .659,  .471}233.25 & \cellcolor[rgb]{ .988,  .686,  .475}229.03 & \cellcolor[rgb]{ .988,  .694,  .475}227.77 & \cellcolor[rgb]{ .988,  .698,  .475}226.79 \\
			  & 1     & \cellcolor[rgb]{ .98,  .545,  .447}252.09 & \cellcolor[rgb]{ .988,  .686,  .475}228.98 & \cellcolor[rgb]{ .992,  .757,  .486}217.41 & \cellcolor[rgb]{ .996,  .796,  .494}210.65 & \cellcolor[rgb]{ .996,  .82,  .498}206.68 & \cellcolor[rgb]{ .996,  .827,  .502}205.26 & \cellcolor[rgb]{ .996,  .835,  .502}204.38 \\
			  & 2     & \cellcolor[rgb]{ .988,  .635,  .463}237.49 & \cellcolor[rgb]{ .996,  .776,  .49}213.70 & \cellcolor[rgb]{ 1,  .851,  .506}201.61 & \cellcolor[rgb]{ 1,  .89,  .514}194.97 & \cellcolor[rgb]{ 1,  .91,  .518}191.94 & \cellcolor[rgb]{ 1,  .922,  .518}190.13 & \cellcolor[rgb]{ 1,  .922,  .518}189.70 \\
			  & 3     & \cellcolor[rgb]{ .992,  .733,  .482}221.08 & \cellcolor[rgb]{ 1,  .878,  .51}197.13 & \cellcolor[rgb]{ .945,  .906,  .514}185.06 & \cellcolor[rgb]{ .871,  .882,  .51}178.38 & \cellcolor[rgb]{ .835,  .875,  .506}175.44 & \cellcolor[rgb]{ .82,  .867,  .506}173.98 & \cellcolor[rgb]{ .808,  .863,  .506}172.86 \\
			  & 4     & \cellcolor[rgb]{ .996,  .82,  .498}206.85 & \cellcolor[rgb]{ .925,  .898,  .51}183.24 & \cellcolor[rgb]{ .796,  .863,  .506}171.99 & \cellcolor[rgb]{ .714,  .839,  .498}164.59 & \cellcolor[rgb]{ .682,  .827,  .498}162.04 & \cellcolor[rgb]{ .663,  .824,  .498}160.08 & \cellcolor[rgb]{ .659,  .824,  .498}159.89 \\
			  & 5     & \cellcolor[rgb]{ 1,  .91,  .518}192.05 & \cellcolor[rgb]{ .745,  .847,  .502}167.42 & \cellcolor[rgb]{ .616,  .808,  .494}156.15 & \cellcolor[rgb]{ .537,  .788,  .49}148.99 & \cellcolor[rgb]{ .51,  .78,  .486}146.79 & \cellcolor[rgb]{ .49,  .773,  .486}145.02 & \cellcolor[rgb]{ .482,  .773,  .486}144.37 \\
			  & 6     & \cellcolor[rgb]{ .918,  .898,  .51}182.78 & \cellcolor[rgb]{ .655,  .82,  .494}159.29 & \cellcolor[rgb]{ .518,  .78,  .486}147.44 & \cellcolor[rgb]{ .443,  .761,  .482}140.73 & \cellcolor[rgb]{ .408,  .749,  .482}137.81 & \cellcolor[rgb]{ .392,  .745,  .482}136.23 & \cellcolor[rgb]{ .388,  .745,  .482}135.79 \\
		\end{tabular}
	
	 \caption{The average number of hours that stations are either full or empty, for each combination of morning and evening part-time shifts.}
	\label{tab:parttime}
 \end{table}

Various conclusions can be drawn from Figure \ref{fig:comparisonparttime}. In particular, anticipating the starting time of the earliest shift from 7:00 to 6:00 allows us to better deal with the earliest part of the morning shortage. This, together with the higher number of vehicles performing the relocation, results in a considerable reduction of the morning peak.

\subsection{Experiments on the composition of the fleet of vehicles} \label{sec:fleetv}

\begin{figure}[h]
	\centering
	\includegraphics[width=0.45\linewidth]{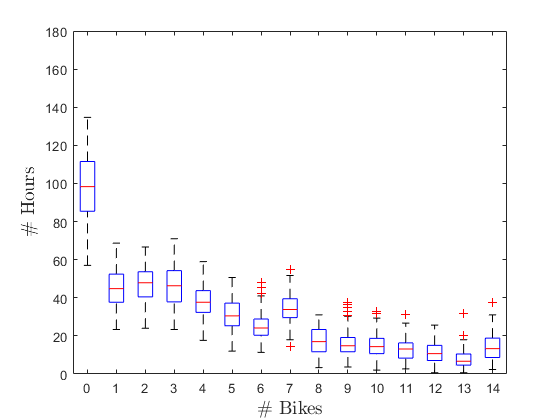}
	\includegraphics[width=0.45\linewidth]{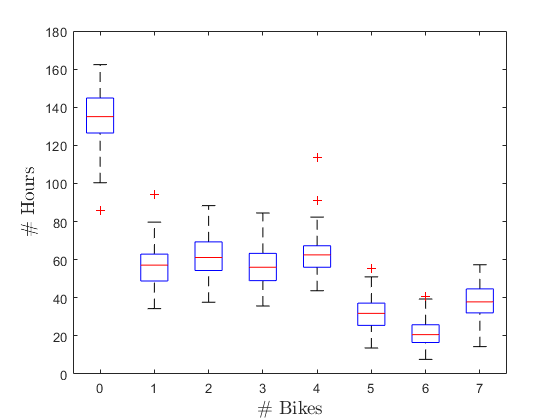}
	\caption{Distribution of the load of the vehicles for the current vehicle fleet and shift configuration (left) and with vehicles with capacity of 7 bikes (right).}
	\label{fig:distload}
\end{figure}

Observing the number of bikes loaded in the vehicles during the day can help to size the fleet of the vehicles. Figure \ref{fig:distload} (left) reports the distribution of the load of the vehicles for the current vehicle fleet and shift configuration: on the x-axis the number of bikes is reported while on the y-axis the box-plots of the distributions of the time that vehicles are loaded with that amount of time across the scenarios are reported. In this configuration, the average load of the vehicles is 4.32 bikes and, as reported in Section \ref{sec:cs}, the average number of hours that stations are either empty or full is 176.03.

Given the distribution of the load of the vehicles when the capacity is of 14 bikes, the company might want to test the use of smaller vehicles with lower capacity. Testing the same setting changing only the capacity of the vehicles from 14 to 7 results in an increase of the hours that stations are either empty or full station of around 0.8\% to 177.42, with an average load of around 2.5 bikes. This increase, however, is not statistically significant. The distribution of the load of the vehicles is reported in Figure \ref{fig:distload} (right).

It is worth pointing out that smaller vehicles are perhaps more nimble and therefore might result in faster travel times, parking times, and allow for faster loading and unloading operations. On the other side, larger vehicles could provide higher resilience to exceptional circumstances. These effects have not been considered in the presented results.

The results presented in this section match the empirical experiments and measurements carried out by the company with respect to the ability of smaller vehicles to fulfill the workload of the current fleet of the company while allowing for more agile operations. These results also support the transition to a fleet of electric vehicles.

\subsection{Experiments on the sizing of the stations} \label{sec:sizes}

\begin{figure}[h]
	\centering
	\includegraphics[width=0.7\linewidth]{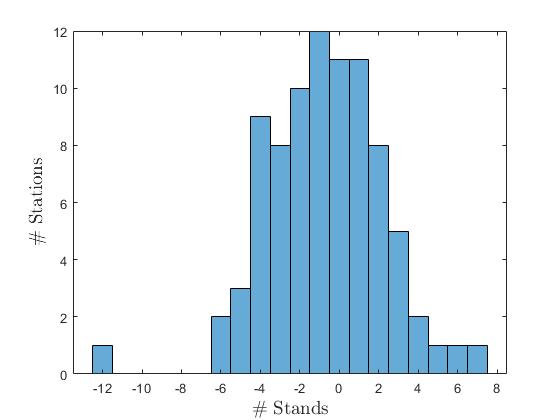}
	\caption{The histograms of the changes in the number of stands.}
	\label{fig:betterbikestands}
\end{figure}

Observing the evolution of the number of bikes in each station can help in sizing the number of stands each station is equipped with. It is not rare that stands are moved from one station to another, or that new stands are bought and installed in the system. To compute the new value for $u_i$ a simple rule could be considered: $ u_i = \max(\hat{f}_{i}(0,H)) - \min(\hat{f}_{i}(0,H)) + 2 \delta$, with $\delta \in \mathbb{N}$. This sets the number of stands available in each station equal to the span of the balance of the requests during the day plus a margin $\delta$ to account for the variability of the requests both on the lower and upper bound. An example of the introduction of this changes is shown in Figure \ref{fig:betterbikestands}, with $\delta = 1$. In this case, 63 stands are to be removed from the stations (68 stands are added and 131 are removed from and to stations). One station appears to be considerably oversized. This station is close to the emergency room and the medicine faculty. As for any practical decision, further considerations should be made before actually implementing the changes presented in this section.

To properly frame the benefits of the sizing of the stations it is worth pointing out that the main cost associated with the stands is their acquisition and, to a lesser degree, their installation. Furthermore, it is also worth pointing out that an undersized station has a greater chance of requiring the users that wish to return a bike to use a lock. This is a cause of discomfort for the users (operators) that have to manually lock (unlock) the bikes.
For these reasons, we test the case where the number of removed stands in each station is halved (rounded to the ceil in the case of odd numbers) with respect to the values computed with the rule reported at the beginning of this section. This results in a net of 14 stands added to the system, an increase of around 1.7\% from the current 827 stands. This increase is in line with the short-term plan of the company to expand and improve the system.
Implementing the changes resulting from this calculation results in 259.81 hours that stations are either full or empty, an improvement of about 5.5\% from the current status of the system when no relocation is performed. In the case of the current shifts of the operators, the result is of 173.41 hours, an improvement of around 1.5\%.

\subsection{Experiments on the sizing of the fleet of bikes} \label{sec:bikefleet}

\begin{figure}[h]
	\centering
	\includegraphics[width=0.95\linewidth]{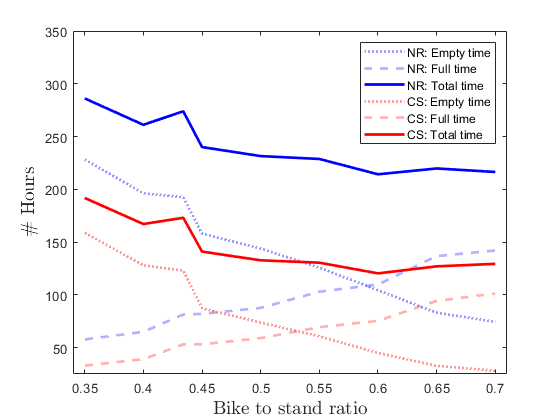}
	\caption{Results for different bikes to stands ratio.}
	\label{sec:bikeratio}
\end{figure}

As reported at the beginning of Section \ref{sec:casestudy}, around 360 bikes are typically available in the system, while 827 stands are available overall. This results in a bike-to-stand ratio of 0.435. As it can be observed from the result provided in Section \ref{sec:numopt}, the number of hours in which stations are empty greatly exceeds the one where stations are full. To test whether it would be beneficial to increase the number of bikes in the system, tests are presented with a decreasing and increasing bike-to-stand ratio of values in the interval $[0.35, 0.7]$, with steps of $0.05$. Results are reported in Figure \ref{sec:bikeratio} for different values of the ratio both for the case where no relocation is performed (indicated as NR in the figure) and the one where the shifts of the vehicles are performed according to the current status of the system (indicated as CS in the figure). Ratios lower than $0.434$ have been achieved by removing bikes proportionally to the stock at the beginning of the day. Similarly, ratios higher than $0.434$ have been achieved by adding bikes in stations proportionally to their number of empty stands at the beginning of the day.
Various conclusions can be drawn from Figure \ref{sec:bikeratio}. First, the current distribution of the stock of the bikes appears to be outlying with respect to the result simulated for different stocks of the stations, i.e., all curves reported in Figure \ref{sec:bikeratio} show a "bump" when the ratio is equal to $0.434$. This is due to the fact that the way bikes have been removed or added to the system provides a little advantage with respect to the current situation, adding or removing bikes where it is intuitively convenient to do so. Second, the benefit of bike relocation appears to be consistent for different sizing of the fleet of bikes, as the distance between matching curves, i.e., the two solid, dashed, or dotted lines, is almost constant. Finally, increasing the number of bikes appears to have a greater impact on the number of empty hours than on that of full hours, i.e., the rate at which the number of empty hours decreases is greater than that at which the number of full hours increased.
Based on the considerations presented in this section, the company is currently testing the performance of the system with a greater bike-to-stand ratio. In particular, preliminary tests are focused in the $[0.50 - 0.55]$ interval, with an increase of around 25\% in the number of bikes in the system.

\section{Conclusions} \label{sec:conclusions}

A simulation framework for a station-based BSS has been introduced. The framework provides a way to model the layout of a BSS, i.e., its stations, as well as the fleet of vehicles and the shifts of the operators performing the relocation of the bikes. The requests of a bike to rent or a stand to return a bike are modeled through stochastic processes.
A dynamic bike relocation algorithm is considered within the framework. The algorithm aims at minimizing the expected shortages in stands and bikes at each station considering the forecast of user demand for rentals and returns. When planning the relocation of bikes in the system the entire fleet of vehicles is considered, as opposed to stations being clustered and vehicles being assigned to the relocation of bikes only for one cluster.
The simulation framework originates from a collaboration with Brescia Mobilità, operating the Bicimia BSS of the city of Brescia, Italy. Computational results have been presented for the case of Bicimia. These results show that the framework provides valuable insights for the design and improvement of a BSS at a strategic and tactical level. \red{We have reported results discussing the number of vehicles and operators to perform the relocation and their shifts, and the composition of the fleet of vehicles, that is, their capacity. We have presented experiments discussing the sizing of the stations as well as the bike fleet.
In particular, in the tested setting, we have shown that the company is currently allocating the best number of operators in the morning and afternoon shifts. Results on the current shift schedule have shown however that eliminating the current overlap between the two shifts would greatly improve the quality of the service and that shorter shifts covering the morning and evening peak would not bring a significant improvement over the current ones. Results on the capacity of the vehicles have shown that the company could deploy smaller vehicles with little change in the current quality of service. A redistribution of the stands of the stations and the addition of few stands in key stations was shown to bear a small improvement to the quality of the system. Finally, we have shown that the bike to stand ratio could be increased to reduce the number of hours that stations are either empty or full.}

Several research directions remain to be studied. Different relocating policies for the vehicles could be tested, to assess how decisions at the operational level impact decisions at the strategic and tactical ones. Another interesting direction is the expansion of the horizon considered by the framework to multiple days. This would allow us to test how decisions regarding the relocation of bikes in the current day impact the performance of the system in the future. Finally, a real-time tool for the dynamic routing of the relocating vehicles would be of relevant practical interest.

\bibliographystyle{abbrvnat}
\bibliography{biblio}

\end{document}